\documentclass[12pt]{article}
\usepackage{amsfonts}
\usepackage{amsmath}
\usepackage{color}

\newtheorem{theorem}{Theorem}

\newenvironment{proof}[1][Proof]{\noindent\textbf{#1.} }{\ \rule{0.5em}{0.5em}}
\input{tcilatex}
\begin{document}

\title{A note on the well posed anisotropic discrete BVP's}
\author{Marek Galewski}
\maketitle

\begin{abstract}
Using the direct method of the calculus of variations we investigate the
existence, uniqueness and continuous dependence on parameters for solutions
of second order discrete anisotropic equations with Dirichlet boundary
conditions.
\end{abstract}

\section{Introduction}

Since difference equations serve as mathematical models in diverse areas,
such as economy, biology, physics, mechanics, computer science, finance -
see \cite{elyadi}- it is of interest to know the conditions which guarantee
the existence and uniqueness of solutions and their dependence on
parameters. Problems satisfying all three conditions are called well-posed.\
We consider an anisotropic difference equation 
\begin{equation}
\left\{ 
\begin{array}{l}
-\Delta \left( h\left( k-1\right) |\Delta x(k-1)|^{p(k-1)-2}\Delta
x(k-1)\right) =\lambda f(k,x(k),u\left( k\right) ),\bigskip \\ 
x(0)=x(T+1)=0,%
\end{array}%
\right.  \label{zad}
\end{equation}%
where $\lambda >0$ is a numerical parameter, $f:Z[0,T+1]\times \mathbb{R}%
\times \mathbb{R}\rightarrow \mathbb{R}$, $Z[a,b]$ for $a<b$, $a,b$
integers, is a discrete interval $\{a,a+1,...,b\};$ $\Delta x(k-1)=x\left(
k\right) -x(k-1)$ is the forward difference operator; $p,$ $h:Z\left[ 0,T+1%
\right] \rightarrow \mathbb{R}_{+}$, $p^{-}=\min_{k\in Z\left[ 0,T+1\right]
}p\left( k\right) >1$; $h^{-}>0$; $u:Z\left[ 1,T\right] \rightarrow \mathbb{R%
}$ is a parameter. A solution to (\ref{zad}) is a function $%
x:Z[0,T+1]\rightarrow \mathbb{R}$\ which satisfies the given equation and
the associated boundary conditions. Solutions are investigated in a space $%
H=\{x:Z[0,T+1]\rightarrow \mathbb{R}:x(0)=x(T+1)=0\}$ which is a Hilbert
space when considered with a norm $||x||=\left( \sum_{k=1}^{T+1}|\Delta
x(k-1)|^{2}\right) ^{1/2}$. With fixed function $u$ and $\lambda >0$ the
functional corresponding to (\ref{zad}) is 
\begin{equation*}
J_{u}(x)=\sum_{k=1}^{T+1}\frac{h\left( k-1\right) }{p(k-1)}|\Delta
x(k-1)|^{p(k-1)}-\lambda \sum_{k=1}^{T}F(k,x\left( k\right) ,u(k)),
\end{equation*}%
where $F(k,x,u)=\int_{0}^{x}f(k,t,u)dt$. $J_{u}:H\rightarrow \mathbb{R}$ is
continuously G\^{a}teaux differentiable and equating its G\^{a}teaux
derivative $J_{u}^{^{\prime }}$\ to zero%
\begin{equation*}
\begin{array}{l}
\left\langle J_{u}^{^{\prime }}(x),y\right\rangle =\sum_{k=1}^{T+1}h\left(
k-1\right) |\Delta x(k-1)|^{p(k-1)-2}\Delta x(k-1)\Delta y(k-1)\bigskip \\ 
-\lambda \sum_{k=1}^{T}f(k,x(k),u\left( k\right) )y\left( k\right) =0\text{
for }y\in H%
\end{array}%
\end{equation*}%
provides a weak solutions to (\ref{zad}). Summing by parts we see that a
weak solution is a strong one- compare with \cite{uni1}, \cite{KoneOuro},
where the weak solutions are obtained. The uniqueness of a solution is
implied by the uniqueness of a critical point and this is turn is guaranteed
by strict convexity. The assumptions leading to the existence and uniqueness
suffice to prove the continuous dependence on parameters.\smallskip

Continuous versions of (\ref{zad}) are known to be mathematical models of
various phenomena arising in the study of elastic mechanics (see \cite{B})
and electrorheological fluids (see \cite{A}). While the research of
continuos anisotropic problems has been very abundant (see \cite{hasto}),
the investigations within their discrete counterparts have only begun \cite%
{uni1}, \cite{KoneOuro}, \cite{MRT} \ - where some tools from the critical
point theory are applied in order to get the existence and multiplicity of
solutions. None of theses sources considers the well-posed problems. In \cite%
{uni1} some uniqueness results are given, but these are provided for some
special type of nonlinear terms and for weak solutions. We have already
undertaken investigations concerning dependence on parameters for discrete
problems in \cite{galewskiCOERCIVE}. Some uniqueness results for discrete
problems can be found in \cite{uni2}, where different approach is applied.

\section{Auxiliary results}

For the following estimations see \cite{MRT}. There exist $C_{1},C_{2}>0$
such that, 
\begin{equation}
\sum_{k=1}^{T+1}|\Delta x(k-1)|^{p(k-1)}\geq C_{1}||x||^{p^{-}}-C_{2}\text{
for every }x\in H\text{ with }||x||\geq 1.  \label{(eq0)}
\end{equation}%
There exists $c_{m}>0$ such that for any $x\in H$ and any $m\geq 2$%
\begin{equation}
\sum_{k=1}^{T}|x(k)|^{m}\leq c_{m}\sum_{k=1}^{T+1}|\Delta x(k-1)|^{m}\text{
and}  \label{eq1}
\end{equation}%
\begin{equation}
(T+1)^{\frac{2-m}{2m}}||x||\leq \left( \sum_{k=1}^{T+1}|\Delta
x(k-1)|^{m}\right) ^{1/m}\leq (T+1)^{\frac{1}{m}}||x||.  \label{relation_m}
\end{equation}

\begin{theorem}
\cite{Ma}\label{LematCritPoint}Let $E$ be a reflexive Banach space. Let $%
J:E\rightarrow \mathbb{R}$, $J\in C^{1}\left( E,\mathbb{R}\right) $, be
weakly lower semi-continuous, coercive and strictly convex. Then there
exists a unique point $x_{0}\in E$ such that $\inf_{x\in E}J(x)=J(x_{0})$
and $J^{^{\prime }}(x_{0})=0$.
\end{theorem}

\section{Main Result}

In this note we assume that

\textbf{H1 }\textit{there exist\ functions} $a:Z\left[ 1,T\right]
\rightarrow \mathbb{R}_{+}$\textit{, }$b:Z\left[ 1,T\right] \rightarrow 
\mathbb{R}$, $q:Z\left[ 1,T\right] \rightarrow \left( 1,+\infty \right) $%
\textit{\ such that }%
\begin{equation*}
|f(k,x,u)|\leq a(k)|x|^{q(k)}+b(k)\text{ \textit{for all} }x,u\in \mathbb{R}%
\text{ \textit{and all} }k\in \left[ 1,T\right] \text{;}
\end{equation*}%
\ \ \qquad \newline

\textbf{H2}\textit{\ for any fixed }$k\in Z\left[ 1,T\right] $ \textit{%
function} $x\rightarrow f\left( k,x,u\right) $\textit{\ is nonincreasing for
all }$u\in \mathbb{R}$.\smallskip

\textbf{H3} $f\left( k,0,u\right) \neq 0$\textit{\ for at least one }$k\in Z%
\left[ 1,T\right] $\textit{\ and for all }$u\in \mathbb{R}$\textit{.}

\begin{theorem}
\label{t37} Assume that conditions \textbf{H1}-\textbf{H3} hold with either $%
p^{-}>q^{+}+1$ and $\lambda >0$ or $p^{-}=q^{+}+1$ and $\lambda <\frac{%
C_{1}h^{-}\left( q^{-}+1\right) }{p^{+}a^{+}c_{q^{+}+1}}$. Then for each $u:Z%
\left[ 1,T\right] \rightarrow \mathbb{R}$ problem (\ref{zad}) has exactly
one nontrivial solution $x_{u}$. Let $\left\{ u_{n}\right\} _{n=1}^{\infty }$
be a convergent sequence of parameters, where $\lim_{n\rightarrow \infty
}u_{n}\left( k\right) =\overline{u}\left( k\right) $ for $k\in Z\left[ 1,T%
\right] $. For a sequence $\left\{ x_{u_{n}}\right\} _{n=1}^{\infty }$ of
solutions to problem (\ref{zad}) corresponding to $u_{n}$, there exists a
convergent subsequence $\left\{ x_{u_{n_{i}}}\right\} _{i=1}^{\infty }$ such
that its limit $\overline{x}$ solves (\ref{zad}) for $u=\overline{u}$.
\end{theorem}

\proof%
For the existence and uniqueness we apply Theorem \ref{LematCritPoint}. Let $%
u$ be fixed. Assumption \textbf{H3} guarantees that all solutions must be
nontrivial. $J_{u}\in C^{1}\left( H,\mathbb{R}\right) $ and it is strictly
convex since by \textbf{H2} the nonlinear terms are convex and since the
terms connected with the difference operator are strictly convex. In order
to show the coercivity take $||x||\geq 1$ with $||x||_{C}=\max_{k\in Z\left[
1,T\right] }\left\vert x\left( k\right) \right\vert \geq 1$ and observe
using (\ref{eq1}), (\ref{relation_m}) that%
\begin{equation*}
\begin{array}{l}
\left\vert \sum_{k=1}^{T}F(k,x(k),u\left( k\right) )\right\vert \leq
\sum_{k=1}^{T}\left( \frac{a(k)|x(k)|^{q(k)+1}}{q(k)+1}+b(k)x(k)\right) \leq
\bigskip \\ 
\sum_{k=1}^{T}\left( \frac{a^{+}}{q^{-}+1}|x(k)|^{q^{+}+1}+b^{+}|x(k)|%
\right) \leq \frac{a^{+}c_{q^{+}+1}}{q^{-}+1}%
(T+1)||x||^{q^{+}+1}+b^{+}c_{1}(T+1)||x||.%
\end{array}%
\end{equation*}%
Further by (\ref{(eq0)}) we get 
\begin{equation}
J_{u}(x)\geq \frac{C_{1}h^{-}}{p^{+}}||x||^{p^{-}}-\lambda \frac{%
a^{+}c_{q^{+}+1}}{q^{-}+1}(T+1)||x||^{q^{+}+1}-\lambda
b^{+}c_{1}(T+1)||x||-C_{2}.  \label{est_coer}
\end{equation}%
Thus $J_{u}(x)\rightarrow \infty $ as $||x||\rightarrow \infty $ in case $%
p^{-}>q^{+}+1$ and also $J_{u}(x)\rightarrow \infty $ as $||x||\rightarrow
\infty $ in case $p^{-}=q^{+}+1$ since $\lambda <\frac{C_{1}h^{-}\left(
q^{-}+1\right) }{p^{+}a^{+}c_{q^{+}+1}}$. So problem (\ref{zad}) has exactly
one, nontrivial solution $x_{u}\in H$. \bigskip

Let $\left\{ x_{u_{n}}\right\} _{n=1}^{\infty }$ be the sequence of
solutions corresponding to $\left\{ u_{n}\right\} _{n=1}^{\infty }$. Suppose
that $\left\{ x_{u_{n}}\right\} _{n=1}^{\infty }$ is unbounded. Then $%
||x_{u_{n}}||\rightarrow \infty $ as $n\rightarrow \infty $. Note also that%
\begin{equation}
J_{u}(x_{u_{n}})=\inf_{x\in H}J_{u_{n}}(x)\leq J_{u_{n}}(0)=0.  \label{Jmn}
\end{equation}%
Since $\left\{ x_{u_{n}}\right\} _{n=1}^{\infty }$ is assumed unbounded
there exists $N_{0}$ such that $||x_{u_{n}}||\geq 1$ and $||x||_{C}\geq 1$
for $n\geq N_{0}$. Now by (\ref{est_coer}) together with (\ref{Jmn}) we see
that for $n\geq N_{0}$ 
\begin{equation*}
\frac{C_{1}h^{-}}{p^{+}}||x_{u_{n}}||^{p^{-}}-\lambda \frac{a^{+}c_{q^{+}+1}%
}{q^{-}+1}(T+1)||x_{u_{n}}||^{q^{+}+1}-\lambda
b^{+}c_{1}(T+1)||x_{u_{n}}||\leq C_{2}.
\end{equation*}%
Thus $\left\{ x_{u_{n}}\right\} _{n=1}^{\infty }$ must be bounded and we
reach a contradiction. Hence there exists a constant $\gamma >0$ such that $%
||x_{u_{n}}||\leq \gamma $ for $n\in 
\mathbb{N}
$ and there exists a convergent subsequence $\left\{ x_{u_{n_{i}}}\right\}
_{i=1}^{\infty }$ whose limit we denote by $\overline{x}$. Let $\left\{
u_{n_{i}}\right\} _{i=1}^{\infty }$ be the corresponding subsequence of
parameters which obviously converges to $\overline{u}$. For these
subsequences we have for$\ k\in Z[1,T]$%
\begin{equation*}
\begin{array}{l}
-\Delta \left( h\left( k\right) |\Delta x_{u_{n_{i}}}(k-1)|^{p(k-1)-2}\Delta
x_{u_{n_{i}}}(k-1)\right) =f\left( k,x_{u_{n_{i}}}(k),u_{n_{i}}(k)\right) ,%
\text{ }\bigskip \\ 
x_{u_{n_{i}}}(0)=x_{u_{n_{i}}}(T+1)=0.%
\end{array}%
\text{ }
\end{equation*}%
Taking limits to both sides of the above relation, we see by continuity that
(\ref{zad}) holds with $x=\overline{x}$ and $u=\overline{u}$.%
\endproof%

\section{Conclusions}

Other variational discrete boundary value problems can be tackled by our
approach provided the assumptions imposed allow for the direct method to be
applied. We also can double easily our results by using strict concavity and
anti-corecivity of a functional $J_{u}^{1}$ 
\begin{equation*}
J_{u}^{1}(x)=\lambda \sum_{k=1}^{T}F(k,x\left( k\right)
,u(k))-\sum_{k=1}^{T+1}\frac{h\left( k\right) }{p(k-1)}|\Delta
x(k-1)|^{p(k-1)}.
\end{equation*}%
Theorem \ref{t37} remains valid with \textbf{H1, H3} retained, with \textbf{%
H2} replaced by\smallskip

\textbf{H4}\textit{\ for all }$k\in Z\left[ 1,T\right] ,$ $u\in \mathbb{R}$ 
\textit{function} $x\rightarrow f\left( k,x,u\right) $\textit{\ is
nonincreasing}\smallskip \newline
and with the following assumptions upon $p$ and $q$ 
\begin{equation*}
q^{-}+1>p^{+},\text{ }\lambda >0\text{ or }q^{-}+1=p^{+},\text{ }\lambda
>(T+1)^{\frac{1-p^{-}}{2p^{-}}}\frac{h^{-}\left( q^{-}+1\right) }{%
p^{+}a^{+}c_{q^{+}+1}}.
\end{equation*}

\begin{tabular}{l}
Marek Galewski \\ 
Institute of Mathematics, \\ 
Technical University of Lodz, \\ 
Wolczanska 215, 90-924 Lodz, Poland, \\ 
marek.galewski@p.lodz.pl%
\end{tabular}

\end{document}